\documentclass[11pt]{article}
\usepackage{latexsym}
\usepackage[dvips]{epsfig}

\newcommand{\qed}{\mbox{$\Diamond$}\vspace{\baselineskip}}

\newtheorem{theorem}{Theorem}[section]
\newtheorem{proposition}[theorem]{Proposition}
\newtheorem{lemma}[theorem]{Lemma}
\newtheorem{definition}[theorem]{Definition}
\newtheorem{corollary}[theorem]{Corollary}
\newtheorem{example}[theorem]{Example}

\newcommand{\vanish}[1]{}
\begin{document}

\title{Comparing algorithms for sorting with $t$ stacks in series}
\author{Rebecca Smith \thanks{
        Department of Mathematics,
        University of Florida,
        Gainesville FL 32611-8105.}}
\maketitle

\date{}
\begin{abstract}
We show that the left-greedy algorithm is a better algorithm than the 
right-greedy algorithm for sorting permutations using $t$ stacks in 
series when $t>1$.  We 
also supply a method for constructing some permutations that can be sorted 
by $t$ stacks in series and from this get a lower bound on the number of 
permutations of length $n$ that are sortable by $t$ stacks in series. 
 Finally 
we show that the left-greedy algorithm is neither optimal nor defines a closed
class of permutations for 
$t>2$.
\end{abstract}

\section{Introduction}

The ideas of stack sorting and sorting by stacks in series were introduced 
by Knuth in \cite{4}.  In this he showed that the permutations of length 
$n$ which could be sorted by one stack is the $n^{th}$ Catalan number.  
These permutations are the ones that avoid the permutation 231.

\begin{definition} We say that a permutation is {\em sortable by a single 
stack} 
if starting at the beginning of a permutation we can put the entries in a 
stack and remove them from the stack to the output as necessary by a 
``last in, first out'' method to get the identity permutation while 
always keeping the elements in the stack in an increasing order from 
top to bottom.
\end{definition}

\begin{example} {\em  The permutation 4132 is stack sortable.  

\begin{figure}[ht]
 \begin{center}
  \epsfig{file=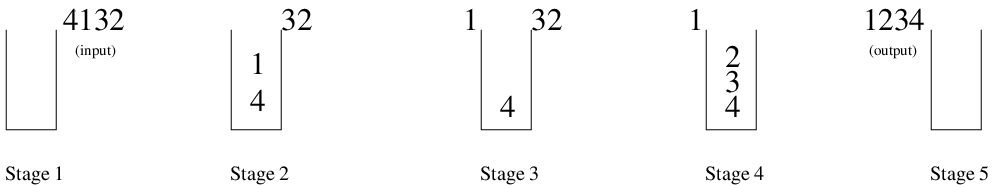}
  \caption{Sorting 4132}
  \label{4132}
 \end{center}
\end{figure}
}
\end{example}

So when we refer to {\em stack sorting using $t$ stacks 
in series}, we use the same monotonicity rules for each stack, but once 
the elements leave a stack, they go directly into the next stack.  
The elements are put into the output when they have passed through 
all stacks.  (The elements pass through the stacks from right to left.)

In \cite{5}, West introduced a right-greedy algorithm to sort permutations 
by $t$ stacks in series.  (The right-greedy algorithm always takes the 
rightmost legal move.)  He showed that the permutations that it sorted 
are exactly those that are $t$-stack sortable, that is, the 
permutations that are sorted after running them through one stack $t$ times.  
(When running the permutation through the stack each time, we use the optimal 
algorithm for 1-stack sorting, that is, only remove elements from the stack 
if the next element cannot enter the stack without violating the increasing 
order condition.)

Later, Atkinson, Murphy, and Ru$\check{s}$kuc used a left-greedy algorithm 
(always taking the leftmost legal move) 
for sorting permutations using two stacks in series in \cite{1}.  They 
showed that using this 
algorithm would sort any permutation that could be sorted by 2 stacks in 
series.  This allowed them to count the number of such permutations of a 
given length $n$.  They also determined that these sortable permutations 
were a closed class.  That is, a permutation was sortable by two stacks 
in series (using the left-greedy algorithm) if and only if 
it avoided an infinite basis of unsortable permutations.

For more information on stack sorting, the reader is encouraged to consult 
\cite{2} or \cite{3}.

In this paper, we show any permutation that can be sorted by the 
right-greedy algorithm on $t$ stacks in series can be sorted by the 
left-greedy algorithm on $t$ stacks in series.  A lower bound on the 
number of permutations that can be sorted by $t$ stacks in series is shown 
as well.  We  present examples 
showing that there are many permutations that can be sorted by the 
left-greedy algorithm, but not the right-greedy algorithm on $t$ stacks 
in series for $t>1$.  Also, examples will be presented showing that 
the left-greedy algorithm is not optimal for more than two stacks in 
series and the permutations that are sorted by the left-greedy algorithm 
for more than two stacks are not a closed class.

\section{The Main Result}

The fact that any permutation that can be sorted by the right-greedy 
algorithm will also be sorted by the left-greedy algorithm on $t$ stacks 
in series will follow from the subsequent lemma.  It is already known 
for the case when $t=1$, since both algorithms are the same and when 
$t=2$ since the left-greedy algorithm is again known to be optimal.

\textbf{Note:}  We will say that an algorithm {\em fails} when we cannot 
move any element to a stack without violating the increasing condition on the 
stacks, nor can we move the next element of the identity permutation into the 
output.

\begin{definition}  Define the
 {\em $i^{th}$ critical moment} 
to be the time right before the element $p_i$ would enter the stacks or
 when the algorithm fails, whichever happens first.

\end{definition}

\begin{lemma}

Suppose we have $t$ stacks in series and let $p=p_1p_2...p_n$ be a
 permutation.
For every $i=1,2,...n$, at the $i^{th}$ critical moment each element of 
the permutation
is at least as far left by the left-greedy algorithm as it 
would be by the right-greedy algorithm at its $i^{th}$ critical moment.
\newline
\noindent
\textbf{Note:}  The output may be considered the furthest left position 
and the input would be the furthest right position.

\end{lemma}

\begin{figure}[ht]
 \begin{center}
  \epsfig{file=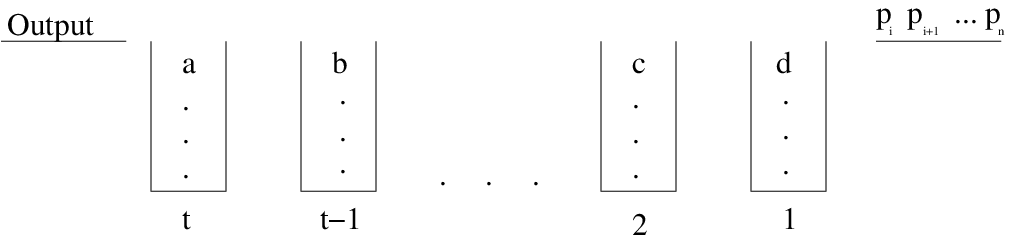}
  \caption{Stacks using the right-greedy algorithm immediately before 
$p_i$ enters}
  \label{rightstacks1}
 \end{center}
\end{figure}

\noindent 
\textbf{Proof:} We use induction on the $i^{th}$ critical moment.
It is clear that the lemma is true for when $i=1$ and when $i=2$.  
So suppose it is true for $i$ and show for $i+1$.
\newline
\indent
If $p_i$ cannot enter the stacks by the left-greedy algorithm, then let $p_j$ 
be the first element of the permutation that could not enter the stacks.  If 
the left-greedy algorithm failed before $p_j$ entered the stacks for $j<i$, 
then the right greedy algorithm must have failed before $p_j$ could enter 
also.  Otherwise the induction hypothesis condition would be false since 
$p_j$ would be further left at the $i^{th}$ critical moment.  If $j=i$, 
then there 
exists a $\gamma < p_i$ in the rightmost stack when using the left-greedy 
algorithm.  However if the right-greedy 
algorithm got even this far without being stuck, then $\gamma$ must be in 
the rightmost stack by this algorithm as well by the induction hypothesis.  
Therefore the right-greedy algorithm fails by this point also, so the 
lemma trivially holds for all subsequent $p_i$'s including $p_{i+1}$.
\newline
\indent
So we may assume that $p_i$ can enter the stacks by the left-greedy 
algorithm.  Also, since we are still trivially done if the right-greedy 
algorithm fails before $p_i$ enters the stacks, we may assume $p_i$ can 
enter the stacks by the right-greedy algorithm as well. 
\newline
\indent
So when using the induction hypothesis, we may assume that the stacks 
using the right-greedy algorithm are as in Figure 
\ref{rightstacks1}, with all elements being at least as far left in the 
stacks using the left-greedy algorithm as those using the right-greedy 
algorithm. 

\begin{figure}[ht]
 \begin{center}
  \epsfig{file=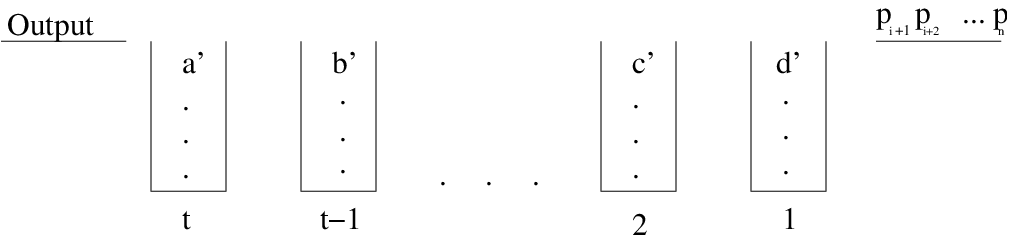}
  \caption{Stacks when the left-greedy algorithm has been carried out 
until the $(i+1)^{st}$ critical moment.}
  \label{leftstacks2}
 \end{center}
\end{figure}

\begin{figure}[ht]
 \begin{center}
  \epsfig{file=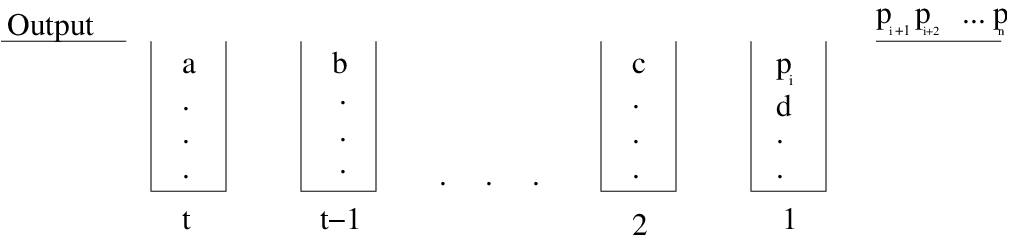}
  \caption{Stacks using the right-greedy algorithm immediately after 
$p_i$ enters.}
  \label{rightstacks2}
 \end{center}
\end{figure}

When $p_i$ enters the first stack by each algorithm, we still have the 
condition that each element is at least as far left by the left-greedy 
algorithm as by the right-greedy algorithm since no elements move left 
except for $p_i$ moving into the first stack in both algorithms.  The 
condition continues to hold if we carry out the left-greedy algorithm until 
it is time for $p_{i+1}$ 
to enter the stacks or the algorithm fails since we are only moving elements 
left in the stacks using the left-greedy algorithm.  So we now have the 
stacks that use the left-greedy algorithm being as shown in Figure 
\ref{leftstacks2} and the stacks that use the right-greedy algorithm as
shown in Figure \ref{rightstacks2}.

\begin{figure}[ht]
 \begin{center}
  \epsfig{file=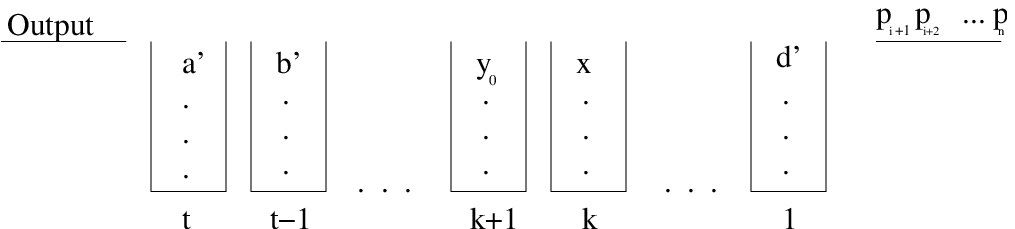}
  \caption{Stacks when the left-greedy algorithm has been carried 
out until the $(i+1)^{st}$ critical moment.}
  \label{leftstacks3}
 \end{center}
\end{figure}

\begin{figure}[ht]
 \begin{center}
  \epsfig{file=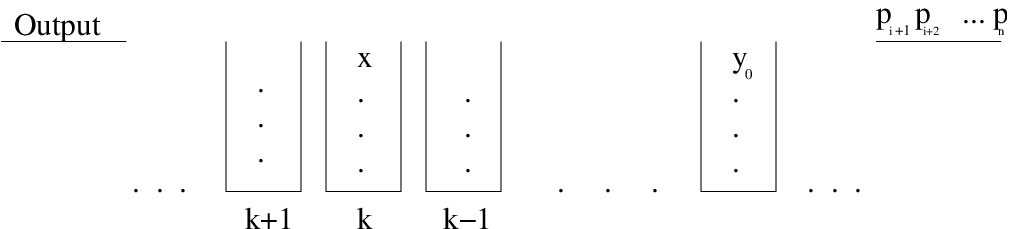}
  \caption{Stacks using the right-greedy algorithm immediately before 
$x$ is to move into the $(k+1)^{st}$ stack.}
  \label{rightstacks3}
 \end{center}
\end{figure}

Now if when we carry out the 
right-greedy algorithm until $p_{i+1}$ should enter or the algorithm 
fails and we end up with all elements no further left than they were by 
the left-greedy algorithm, then we are done.  So suppose there is an element 
$x$ such that $x$ was in the $k^{th}$ stack by the left-greedy 
algorithm at this point, but $x$ moves further left by the right-greedy 
algorithm before $p_{i+1}$ enters.  Choose $x$ to be the 
first element that this happens for (in the order of the steps of the 
right-greedy algorithm) and look at what is happening when $x$ 
moves from the $k^{th}$ stack to the $(k+1)^{st}$ stack.  
\newline
\indent
We first note that if such an element exists, the $(k+1)^{st}$ stack is 
actually a stack and not the output.  This is because
$x$ would not be in the last stack by the left-greedy algorithm 
since any element that has already had all elements smaller than it 
enter the stacks, leaves by the left-greedy algorithm before a new 
element enters.  Thus no element can leave the stacks by the right-greedy 
algorithm that has not also left by the left-greedy algorithm.
\newline
\indent
Since $x$ did not move to the $(k+1)^{st}$ stack by the left-greedy 
algorithm, there must be an element $y_0 <$ $x$ in the 
$(k+1)^{st}$ stack at this point by this algorithm.  As all elements 
are no further left by the right-greedy algorithm at the point when 
we are about to move $x$ from the $k^{th}$ stack to the $(k+1)^{st}$ 
stack, it must be that $y_0$ is to the right of $x$ at 
this point using the right-greedy algorithm.
So the stacks sorted by the left-greedy algorithm must be as in 
Figure \ref{leftstacks3} and the stacks sorted by the right-greedy 
algorithm must be as in Figure \ref{rightstacks3} at this point.

\begin{figure}[ht]
 \begin{center}
  \epsfig{file=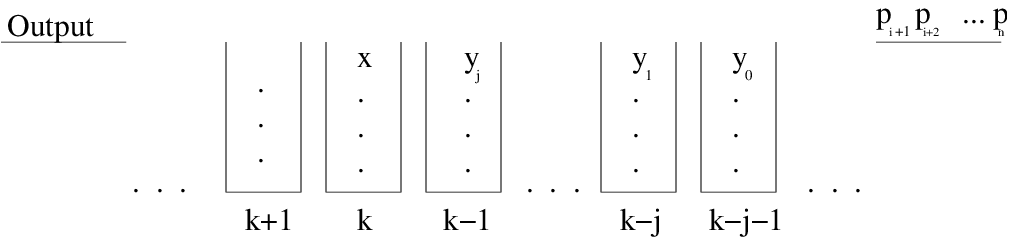}
  \caption{Stacks using the right-greedy algorithm immediately before 
$x$ is to move into the $(k+1)^{st}$ stack.}
  \label{rightstacks4}
 \end{center}
\end{figure}

Because $y_0$ is in the stacks and further right than $x$, 
if $y_0$ could have moved, 
it would have priority over $x$ 
by the right-greedy algorithm.  
If $y_0$ was not the top element of its stack at this point, then the 
smaller element at the top of the stack would have priority over $x$.
So there must be a 
$y_1 < y_0$ at the top of the stack directly to the left of $y_0$.  
Likewise, since $y_1$ is to the right of $x$ and it also 
didn't move, there must be a $y_2 < y_1$ at the top of the stack directly to 
the left of $y_1$ as seen in Figure \ref{rightstacks4}.  This 
process continues until we have $x$ directly to the left of 
$y_j$ for some $j$.  But then moving $y_j$ to the $k^{th}$ 
stack above $x$ would be a further right move than moving 
$x$.  This is a contradiction to the assumption that moving 
$x$ was the next move by the right-greedy algorithm.

So there is no element that moves further left by the right-greedy 
algorithm than by the left-greedy algorithm at the $(i+1)^{st}$ critical 
moment.  Therefore the Lemma is proved.

\qed \newline

\begin{theorem}

Any permutation sorted by the right-greedy algorithm on $t$ stacks in 
series is also sorted by the left-greedy algorithm on $t$ stacks in series.

\end{theorem}

\noindent 
\textbf{Proof:}  This argument follows that of the previous Lemma in 
the case of 
failure at the $i^{th}$ critical moment.  
If the left-greedy algorithm fails, let $p_i$ be the first element of the 
permutation that could not enter the stacks. Then there 
exists a $\gamma < p_i$ in the rightmost stack using the left-greedy 
algorithm.  
However if the right-greedy 
algorithm got even this far without being stuck, then $\gamma$ must be in 
the rightmost stack by this algorithm as well by the previous Lemma.  
Therefore the right-greedy algorithm also fails.

\qed

So to see that the left-greedy algorithm sorts more permutations than the 
right-greedy algorithm when $t>1$, we simply need an example of a 
permutation that is not sortable by t stacks in series using the 
right-greedy algorithm, but is sortable on $t$ stacks using the 
left-greedy algorithm.

\begin{example} {\em For $t>1$, the permutation $(t+1)t(t-1)...32(t+2)1$ is 
not sortable by $t$ stacks in series when using the right-greedy algorithm, 
but is sortable by $t$ stacks in series by the left-greedy algorithm.  In 
fact, the left-greedy algorithm only needs two stacks in series 
to sort this permutation 
regardless of what $t$ is.}
\end{example}

\section{Obtaining some permutations that are sortable by $t$ stacks
 in series}

\begin{lemma} \label{first}

Let $p=p_1p_2...p_{n-1}$ be sortable by the left-greedy algorithm on $(t-1)$ 
stacks in series.
Inserting $n$ anywhere in the permutation will result in a permutation 
$p'$ that is stack sortable by the left-greedy algorithm by $t$ stacks in 
series.

\end{lemma}

\textbf{Proof:}  Since $p$ is sortable on only $t-1$ stacks in series by 
the left-greedy algorithm, we know that if we try to sort our augmented 
permutation $p'$ with the left-greedy algorithm on $t$ stacks, as each 
element in $p'$ is about to enter the stacks until (possibly) after $n$ 
enters, the rightmost stack will be empty.  This is because the 
left-greedy algorithm would always move the elements in the stacks left 
before allowing a new element to enter.  So if we could not get one of 
the original elements of $p$ to move out this first stack (before $n$ 
has the possibility of interfering), $p$ would not be sortable by the 
left-greedy algorithm on $t-1$ stacks.  Because of this, it is clear 
that all elements up to and including $n$ of $p'$ will enter the stacks 
by the left-greedy algorithm.

\vspace{.2 cm}

\emph{Claim:}  If we ignore the moves made by $n$, the left-greedy 
algorithm on $p'$ proceeds exactly the same way as it does on $p$ 
when we try to sort both by $t$ stacks.

\vspace{.2 cm}

\emph{Proof of Claim:}  This is evident until the element $n$ enters the 
stacks when sorting $p'$.  If there is some discrepancy, let moving 
$x$ from the $k^{th}$ stack to the $(k+1)^{st}$ stack be the 
first move made on $p'$ that does not match that on $p$.  Say moving 
$y$ from the $l^{th}$ stack to the $(l+1)^{st}$ stack is the 
corresponding move on $p$.  Because the previous moves all corresponded 
before this point, each element $p_i \neq n$ is in the same position for 
both permutations until we make this move.
\newline
\indent
Note that neither $x$, nor $y$ can be leaving the stacks 
for the output at this time.  This is because both permutations have the 
same output and if one permutation has an element other than $n$ at the 
top of the leftmost stack, the other has the same element there.  So if 
the element should be the next out, it is the leftmost move for both.
\newline
\indent
Case 1:  $k>l$
\newline
Then moving $x$ from the $k^{th}$ stack to the $(k+1)^{st}$ stack 
is the further left move.  But then $x$ should have been able to 
move when the left-greedy algorithm was applied to $p$ as well since 
there can be no additional elements in the $(k+1)^{st}$ stack for $p$ 
that are not there for $p'$.
\newline
\indent
Case 2:  $l>k$
\newline
The same argument as in Case 1 holds because the only extra element that 
could be in the $(l+1)^{st}$ stack for $p'$ is $n$ and since all the 
other elements are smaller, $n$ cannot prevent another element from 
moving left.
\newline
\indent
This completes the proof of the claim.  So since we know $n$ can enter 
the $t$ stacks for the permutation $p'$ using the left-greedy algorithm 
and all elements $1,2,3,...,(n-1)$ can enter the $t$ stacks for the 
permutation $p$ by the left-greedy algorithm, by the claim we know all 
the elements of $p'$ can enter the $t$ stacks when using the left-greedy 
algorithm.  
Hence $p'$ is sortable by $t$ stacks using the left-greedy algorithm.

\qed

It is even easier to show the result is true if we allow any algorithm 
to be used to sort the permutation.

\begin{lemma} \label{second}

Let $p=p_1p_2...p_{n-1}$ be sortable by some algorithm on $(t-1)$ 
stacks in series.
Inserting  $n$ anywhere in the permutation will result in a permutation 
$p'$ that is stack sortable by  $t$ stacks in 
series.

\end{lemma}

\textbf{Proof:}  Since we only needed $t-1$ stacks in series to sort $p$ 
by some algorithm *, we may consider an algorithm on $p'$ that proceeds 
exactly as * does on the elements originally in $p$ except that when 
the elements first come in, they move one stack further left than they 
originally stopped.  That is, for the elements $p_i \neq n$, the first 
stack is essentially skipped.  This leaves the first stack open for $n$ 
to enter.  We leave $n$ in that first stack until all the other elements 
have left the stacks and then finally we move $n$ through the remaining 
$(t-1)$ stacks and into the output, sorting $p'$.  (All other elements 
that enter the stacks after $n$ may easily move past $n$ since 
they are smaller than $n$.)

\qed

This method of constructing permutations that are sortable by $t$ stacks 
in series does not work if we only want to allow the right-greedy algorithm 
to sort the permutation.

\begin{example} {\em Consider the permutation $321$ which is sortable by 
1 stack using the right-greedy algorithm.  Putting $4$ into the third 
position results in the permutation $3241$ which is not sortable by 
the right-greedy algorithm on two stacks in series.}
\end{example}

We can also start with permutations that are sortable by $t$ stacks in 
series to get others that are as well.

\begin{proposition}

Let $p$ be a permutation of length $n-1$ that is sortable by 
$t$ stacks in series using the left-greedy algorithm.  We may insert $n$ 
into any one of the first $t$ positions or into the last position to get 
a new permutation $p'$ that is also sortable by $t$ stacks using the 
left-greedy algorithm.  

\end{proposition}

\textbf{Proof:}  This is clear for when $n$ is put in the last position.  
So suppose $n$ was inserted into one of the first $t$ positions.  It is clear
 that $n$ may enter the stacks at its turn since there are at most $t-1$ 
elements ahead of $n$ and there are $t$ stacks.  The remainder of the 
proof follows exactly as the proof for Lemma 3.1.

\qed

\begin{corollary}
Let $n \geq t$.  Then the number of permutations of length $n$ that are 
sortable by $t$ stacks 
in series (using the left-greedy algorithm) is at least $t!(t+1)^{n-t}$.
\end{corollary}

\textbf{Proof:}
Every permutation of length $t+1$ is sortable by $t$ stacks in series 
using the left-greedy algorithm.  We may insert $t+2$ into any of the first
 $t$ positions or into the last position to get a permutation of length 
$t+2$ that can be sorted using the left-greedy algorithm on $t$ stacks in 
series.  Repeat with $t+3,t+4,...n$.  This 
gives us $(t+1)!(t+1)^{n-t-1} = t!(t+1)^{n-t}$.

\qed

Once again, the previous proposition will also hold for sorting 
permutations on $t$ stacks
 in series using any algorithm that works.  However, we first need the 
following lemma.  

\begin{lemma}
Suppose $p=p_1,p_2,...,p_n$ is a permutation that can be sorted by 
$t$ stacks in series by 
some algorithm.  Then there is an algorithm that never needs to leave an 
empty stack to the left of an element already in the stacks at the $i^{th}$ 
critical moment for any given $i=1,2,...,n$.
\end{lemma}

\textbf{Proof:}  Suppose not.  Then given an algorithm with a minimal 
number of stacks 
left empty at a minimal number of critical moments that does work, 
consider the first critical moment a stack must be left empty.  Let $p_j$ 
be an element in the first nonempty stack to the right of an empty 
stack.  Move $p_j$ one place to the left.  Notice that any algorithm that 
works, can still proceed as before since $p_j$ does not block any element 
that would have been able to move to left of it had it not been moved.  
Hence the algorithm still works with one less empty stack at a critical 
moment.  This is a 
contradiction.

\qed

\begin{proposition}

Let $p$ be a permutation of length $n-1$ that is sortable by 
$t$ stacks in series.  We may insert $n$ 
into any one of the first $t$ positions or into the last position to get 
a new permutation $p'$ that is also sortable by $t$ stacks.  

\end{proposition}

\textbf{Proof:}  This is clearly true when we insert $n$ into the last 
position.  So assume that $n$ was inserted into one of the first $t$ 
positions.  
Since we do not need to leave empty stacks, it is clear that 
we may follow an algorithm that works for $p$ and 
$n$ can enter the stacks as there are at most $t-1$ elements ahead of it.  
After that we may leave $n$ in the first stack until we are done carrying out 
the remainder of the algorithm since $n$ will not interfere with any other 
elements moving into the stacks.  Finally, move $n$ through the remaining 
stacks and into the output to finish sorting $p'$.

\qed

However, we do not always get a permutation that is sortable by the 
right-greedy algorithm when we insert $n$ into one of the first $t$ 
positions of a permutation of length $n-1$ that is sortable by the 
right-greedy algorithm on $t$ stacks in series.  (Inserting an $n$ at the 
end of such a permutation is easily seen to produce another permutation 
that is sortable by $t$ stacks in series using the right-greedy algorithm.)

\begin{example} {\em  The permutation 6372451 is sortable by the
 right-greedy algorithm by three stacks in series, but 68372451 is not.}
\end{example}

\section{Some more information about the left-greedy algorithm for $t >2$}

Although the left-greedy algorithm sorts any permutation that can be 
sorted by $t$ stacks in series for  $t=1$ and $t=2$, the left-greedy 
algorithm is not optimal when $t>2$.

\begin{example} {\em  While it can be seen that the permutation 
$254167...(t+4)3$ is sortable by $t$ stacks in series when $t>2$, 
the left-greedy algorithm fails when $(t+4)$ would need to enter.}
\end{example}

Also, the permutations sorted by the left-greedy algorithm 
are not a closed class when $t>2$.  (Since it is known by \cite{4} that the 
permutations sorted by $t$ stacks in series is a closed class, this is 
actually enough to show that the left-greedy algorithm is not optimal 
for $t>2$.)

\begin{example} {\em  Let $t > 2$.  While it can be seen that the permutation 
$254167...(t+4)3$ is not sortable by $t$ stacks in series using the 
left-greedy algorithm, the permutation $2635178...(t+5)4$ which contains 
$254167...(t+4)3$ is sortable using the left-greedy algorithm by $t$ stacks 
in series.}
\end{example}

\vskip 1 cm 

\begin{center} {\bf Acknowledgment} \end{center}
I wish to thank  Mikl\'os B\'ona for his helpful suggestions 
and encouragement.

\end{document}